



\documentstyle{amsppt}

\topmatter
\title
On the Malcev completion of K\"{a}hler groups
\endtitle
\author
Jaume Amor\'{o}s
\endauthor
\affil
Universitat de Barcelona \\
Dep. d'Algebra i Geometria -
Fac. de Matem\`{a}tiques      \\
Gran Via 585 -
08007 Barcelona, Spain \\
amoros\@cerber.mat.ub.es
\endaffil
\thanks
\flushpar Supported by the Generalitat de Catalunya.
\newline Partially supported by DGCYT grant PB93-0790
\endthanks
\date
29/7/94
\enddate
\endtopmatter

\document

\TagsOnRight
\magnification 1200
\NoBlackBoxes
\baselineskip 13pt

\define\C{\Bbb C}
\define\gdt{\Gamma_2/\Gamma_3}
\define\gr{\text{Gr}\,}
\define\grl{\text{Gr}\,\Cal L_2}
\define\gud{\Gamma_1/\Gamma_2}
\define\gut{\Gamma_1/\Gamma_3}
\define\ld{\Cal L_2}
\define\R{\Bbb R}

\flushpar {\bf Introduction.}

\bigpagebreak

The study of compact K\"{a}hler manifolds made by Hodge and
others shows that a K\"{a}hler structure imposes very strong conditions
on the homotopy type of a compact complex manifold $X$. In particular,
unlike in the case of compact differentiable or closed complex
manifolds, not every finitely presented group $G$ is the fundamental
group of a compact K\"{a}hler manifold. Such groups are called
K\"{a}hler groups.

This note has been inspired by the recent work of F. Johnson and E. Rees
([JR]) and M. Gromov ([G]), showing that free products, and in
particular free groups, are not K\"{a}hler. It has been our purpose to
extend this result and find other restrictions on the
presentations of K\"{a}hler groups. This is done by translating
properties of cup products in $H^*(X)$ into properties of the group
bracket in $\pi_1X$, an idea that came out of [JR], and also by
examining the
Albanese map $X \rightarrow Alb(X)$ after [C]. We describe an algorithm
derived from [St] to compute $\gud G, \gdt G \otimes \R$ from a given
presentation of a group $G$, and use it to give three conditions for the
groups to be K\"{a}hler: The Lie algebra $\grl G$, equivalent to the
holonomy algebra, cannot be free (3.3); one- or two-relator K\"{a}hler
groups either have a torsion abelianized or have a Malcev completion
isomorphic to that of a compact Riemann surface group (4.6); nonfibered
K\"{a}hler groups must satisfy certain lower bounds for the number of
their defining relations, equivalently upper bounds for the rank of
$\gdt G$ (5.7,5.8).

In \S1 we recall the real Malcev completion $G \otimes \R$ of a group
$G$, its equivalent Lie algebra $\Cal L G$, and a 2-step nilpotent Lie
algebra $\grl G \cong (\gud G \otimes \R) \oplus (\gdt G \otimes \R)$,
which is determined by $\gut G/_{\text{Torsion}}$ and is equivalent to
the cup products $H^1(X) \wedge H^1(X) \rightarrow H^2(X)$. This
algebra is actually equivalent to the holonomy algebra of $G$ (cf.
[Ch], [Ko]), and is more convenient for our computations. By
[M2],[DGMS], when $G$ is a K\"{a}hler group the algebra $\grl G$
determines the Malcev completion $G \otimes \R$. In \S2 we
briefly recall Sullivan's 1-minimal model of $X$, its duality with
$\Cal L \pi_1X$, and how the algebra $\grl \pi_1X$ and the product map
$\cup \: H^1(X) \wedge H^1(X) \rightarrow H^2(X)$ are equivalent.

In \S3 we use these results to show that
if $\grl G \cong \grl F_n$, where $F_n$ is a finite rank free group then
$G$ is not K\"{a}hler. This is a
strong quantitative restriction, since the generic group presented
with few relations verifies $\Cal L_2 \cong \grl F_n$ for some
$n$ (see Remark 1.14).

The groups with the simplest presentation after free groups are 1- and
2-relator groups. In Theorem 4.6 we determine the Malcev completions of
1- and 2-relator K\"{a}hler groups, which to a great extent characterize
the groups themselves. The mean to do this is to bound from above $\dim
\gdt \pi_1X \otimes \R$ for any K\"{a}hler group $G$ by a function of
the dimension of the image $Y$ of $X$ by its Albanese map $\alpha \: X
\rightarrow Alb(X)$. A desingularization $\tilde Y$ of $\alpha(X)$ has
been shown by F. Campana ([C]) to verify $\pi_1X \otimes \R \cong
\pi_1 \tilde Y \otimes \R$. It turns out of our work that as $\dim Y$
increases linearly, $\dim H^1(X) \wedge H^1(X)- \dim \gdt \pi_1(X)
\otimes \R$ grows quadratically (Prop. 4.5).

Finally, in \S5, we have established a distincton between fibered and
nonfibered K\"{a}hler groups, and used
the mentioned techniques to give upper bounds for $\dim \gdt G \otimes
\R$ for nonfibered groups, or equivalently lower bounds for their number
of defining relations (Prop. 5.7, Cor. 5.8).

To proof our results and make them effective, we give in \S1 an
algorithm for computing $\gud G \otimes \R$, $\gdt G \otimes \R$ and $
\grl G$ from a given presentation of $G$. This algorithm, which is easy
to
implement by means of the Magnus algebra of free groups, is derived from
a spectral sequence given in [St], and was communicated to the author by
M. Hartl. We use it in Cor. 1.13 and Rmk. 1.14 to show many cases in
which the hypothesis of Thm. 3.3 are fulfilled.

To illustrate our results, we give throughout the paper examples
 of groups that cannot be K\"{a}hler, most of them previously unknown to
the author.

\smallpagebreak

\flushpar {\sl Acknowledgements}: The author wishes to thank V. Navarro
Aznar, for
originally proposing the problem and many subsequent discussions on it.
M. Hartl for his valuable help in \S1, F. Campana and the members
of the Geometry Seminar of the Universitat de Barcelona are also to be
thanked.

\vskip 1cm

\flushpar {\bf \S1. Nilpotent Lie algebras related to a group}

\bigpagebreak

We will recall here the concept of Malcev completion of a group and some
related nilpotent Lie algebras $ \Cal L_n G$. We also give an algorithm
to compute $\grl G$ derived from [St].

\medpagebreak

Let $G$ be a finitely presented group. We can
functorially assign to it a tower of nilpotent groups
$$
\dots \rightarrow \gut G \rightarrow \gud G \rightarrow 1
$$
where $\Gamma_1G=G$, $\Gamma_nG= [\Gamma_{n-1}G,G]$ and
$\Gamma_1/\Gamma_n G = G/ \Gamma_nG$.

A group $G$ is said to be {\it uniquely divisible} when for any pair $g
\in G$, $n \in \Bbb Z$, $g$ has a unique $n$th root in $G$. The category
$n@-\Bbb Q@-\Cal Gr$ of uniquely divisible nilpotent groups is included
in the category $n@-\Cal Gr$ of nilpotent groups, and the inclusion
functor has a left adjoint, the {\bf Malcev completion} functor
$\otimes \Bbb Q \: n@-\Cal Gr \rightarrow n@-\Bbb Q@-\Cal Gr$. The
functor $\otimes \Bbb Q$ is the ordinary tensor product on abelian
groups. For
two alternative ways of defining the Malcev completion, see [HMR]
Part I, or App. A of [Q].

The Baker-Campbell-Hausdorff formula gives a categorical equivalence
between finitely generated groups in $n@-\Bbb Q@-\Cal Gr$ and
finite-dimensional nilpotent $\Bbb Q$-Lie algebras. In the latter
category there are $\otimes \R, \otimes \Bbb C$ functors, and crossing
back and forth in this manner we may define a $\otimes \R$ functor over
$n@-\Cal Gr$. Thus we naturally associate to $G$ a tower
$$
\dots \rightarrow \gut G \otimes \R \rightarrow \gud G \otimes \R
\rightarrow 1
$$
of uniquely divisible nilpotent Lie groups, and its correponding tower
$$
\dots \rightarrow \ld G \rightarrow \Cal L_1 G \rightarrow 0
\tag 1.1
$$
of nilpotent $\R$-Lie algebras.

Denote the lower central series of a Lie
algebra $\Cal L$ as $\Cal L^{(1)}= \Cal L$, $\Cal L^{(n)}= [\Cal
L^{n-1}, \Cal L]$. There is another tower of nilpotent $\R$-Lie
algebras naturally associated to a group $G$: $\text{Gr}_n\,G \otimes
\R= \bigoplus_{i=1}^{n} \Gamma_i/ \Gamma_{i+1} G \otimes \R$, with
bracket induced by the group bracket. We sum up the properties of the
tower of Lie algebras (1.1) that we will apply:

\proclaim{Proposition 1.2}

\item{(i)} The Lie algebras $\Cal L_n G$ have nilpotency class
$\text{nil}\, \Cal L_n G=n$.
\item{(ii)} The tower maps $\Cal L_{n+1} G \rightarrow \Cal L_n G$
induce isomorphisms $\Cal L_{n+1}G/\Cal L_{n+1}G^{(n+1)} \cong \Cal L_n
G$.
\item{(iii)} There are isomorphisms of $\R$-vector spaces $\Cal
L_n^{(n)} G \cong \Gamma_n/\Gamma_{n+1} G \otimes \R$.
\item{(iv)} The graduation of $\Cal L_n G$ by its lower central series
produces a natural tower of isomorphisms of graded Lie algebras
$\text{Gr}\, \Cal L_n G \cong \text{Gr}_n\,G \otimes \R$.
\endproclaim

We will call the Malcev algebra and denote $\Cal LG$ the pronilpotent
algebra $\varprojlim \Cal
L_n G$, which is equivalent to the Malcev completion $G \otimes \R$ by
the Baker-Campbell-Hausdorff formula.

When the group $G$ is the fundamental group of a topological space $X$,
the abelian algebra $\Cal L_1 G$ is just $H_1(X; \R)$. We will
consider in this note the following simplest algebra, $\ld G$, and its
graduate
$\grl G \cong \gud G \otimes \R \oplus \gdt G \otimes \R$, in which
$\ld G$ is included. The algebra $\ld G$ is the quotient of the Malcev
algebra $\Cal L G$ by its third commutator ideal $\Cal L G^{(3)}$, and
is also the quotient of the holonomy algebra of $G$ $\frak g_G$ (cf.
[Ch],[Ko]) by its third commutator ideal.

The groups $G$ we will study will be given by
finite presentations $G= \langle x_1, \dots, x_n \; ; \; r_1, \dots,
r_s \rangle$. This means that $G$ is defined by
$$
1 \longrightarrow N \longrightarrow F \longrightarrow G \longrightarrow
1
\tag 1.3
$$
where $F$ is the free group generated by the generator set $\{ x_1,
\dots, x_n
\}$, and $N$ is the normal subgroup of $F$ spanned by the relation set
$\{ r_1, \dots, r_s \} \subset F$.

We describe the above constructions in a case which is fundamental for
our purposes:

\demo{Example 1.4} Free groups.

Let $G=F_n= F_{\{x_1, \dots, x_n \}}$. Its Malcev completion and Lie
algebras $\gr \Cal L_m F_n$ may be computed by means of its group
algebra (cf. [MKS],[Q],[S]). The conclusion is that, denoting
by $\Cal L(S)$ the free $\R$-Lie algebra spanned by a set $S$, there are
isomorphisms
$$
\gr \Cal L_m F_n \cong \Cal L(\{ X_1, \dots, X_n \})/ \Cal L(\{X_1,
\dots, X_n \})^{(m+1)}
$$
In particular, $\gud F_n \otimes \R \cong \R x_1 \oplus \dots \oplus \R
x_n$, $\gdt F_n \otimes \R \cong \R (x_1,x_2) \oplus \dots \oplus \R
(x_{n-1}, x_n)$, and the brackets in $\grl F_n$ are the group ones in
$\gud F_n$ and zero all others.
\enddemo

The Lie algebra $\grl G$ for a finitely presented $G$ may be obtained
from its presentation and $\grl F$. We will use an algorithm for
computing them derived from [St], where a spectral sequence that
computes all $J_G^m/J_G^{m+1}$ is described, communicated to the author
by M. Hartl.

Consider a group presentation $G= \langle x_1, \dots, x_n \; ; \; r_1,
\dots, r_s \rangle$, which induces the exact sequence given in (1.3).
Let $\R F, \R G$ be the $\R$-group algebras of $F,G$, and denote by
$J_F, J_G$ their respective augmentation ideals. The sequence (1.3)
induces an exact sequence of $\R$-algebras
$$
0 \longrightarrow K \longrightarrow \R F \longrightarrow \R G
\longrightarrow 0
\tag 1.5
$$
where $K$ is the two-sided ideal generated by the $\R$-vector space $D=
\langle r_1-1, \dots, r_s-1 \rangle \subset J_F$. This sequence
restricts to exact sequences
$$
0 \longrightarrow K \longrightarrow J_F^m+K \longrightarrow J_G^m
\longrightarrow 0
\tag 1.6
$$
for all $m \ge 1$. We will compute $J_G/J_G^2$, $J_G^2/J_G^3$ from those
sequences:

\proclaim{Proposition 1.7}
Consider the linear map $f \: \bigoplus_{i=1}^{s} \R r_i
\longrightarrow J_F$ determined by $ r_i \mapsto r_i-1$.
\item{(i)} Let $d_0 \: \oplus \R r_i \rightarrow J_F/J_F^2$ be the
projection of $f$. Then $\text{coker}\, d_0 \cong J_G/J_G^2$.
\item{(ii)} The map $f$ induces a linear map
$$
\align
d_1 \: \ker d_0 \longrightarrow &J_F^2/(J_F^3+ J_F \cdot D+ D \cdot J_F
\\
\sum \lambda_i r_i \longmapsto &\sum \lambda_i (r_i-1)
\endalign
$$
and $\text{coker}\, d_1 \cong J_G^2/J_G^3$.
\endproclaim

\demo{Proof}
(i) The exact sequences of (1.6) induce an isomorphism $J_G/J_G^2 \cong
J_F/ J_F^2+K$. As $K$ is the two-sided ideal spanned by $D$ and $\R F
\cong \R \oplus J_F$, actually $J_F^2+K= J_F^2+D$, and thus $J_G/J_G^2
\cong J_F/J_F^2+D$. By its construction, $\text{Im}\, d_0= D$, and this
proves (i).

(ii) Again by (1.6) we have
$$
J_G^2/J_G^3 \cong \left( J_F^2/(J_F^2 \cap
K) \right) / \left( J_F^3/(J_F^3 \cap K) \right) \cong J_F^2/( J_F^3+
J_F^2 \cap K)
$$
The last denominator is $J_F^3+ J_F^2 \cap K = J_F^3+ J_F \cdot D+ D
\cdot J_F+ D \cap J_F^2$. Obviously $f(\ker d_0) \subset J_F^2$ and thus
$d_1$ is well defined. Moreover, its image is precisely $D \cap J^2_F$,
and (ii) follows from this.
\qed \enddemo

We now relate the computed modules $J_G/J_G^2, J_G^2/J_G^3$ with the
sought ones $\gud, \gdt G \otimes \R$ applying a theorem by D. Quillen
([Q2]):

\proclaim{Theorem}
Let $G$ be a group, $k$ a field of characteristic zero, $k G$ the group
algebra and $j \: \oplus \Gamma_n/\Gamma_{n+1} G \otimes k
\rightarrow \oplus J_G^n/J_G^{n+1}$ given by $g \mapsto g-1$ over the
homogeneous components.

Then $j$ induces an isomorphism of algebras $U(\oplus
\Gamma_n/\Gamma_{n+1} G \otimes \R) \cong \oplus J_G^n/ J_G^{n+1}$.
\endproclaim

In the cases $n=1,2$ this means:

\proclaim{Corollary 1.8}
\item{(i)} $\gud G \otimes \R \cong J_G/J_G^2$.
\item{(ii)} Consider the inclusion $J_G \wedge J_G \hookrightarrow
J_G^2$ given by $x \wedge y \mapsto xy-yx$. Then
$$
\gdt G \otimes \R \cong (J_G \wedge J_G+ J_G^3)/J_G^3 \subset
J_G^2/J_G^3
$$
\endproclaim

Corollary 1.8 allows us to adapt the algorithm of Prop. 1.7 to compute
$\gud, \gdt G \otimes \R$:

\proclaim{Lemma 1.9}
The image of the restriction $f \: \ker d_0 \rightarrow J_F^2$ lies in
$J_F \wedge J_F + J_F^3 \subset J_F^2$.
\endproclaim

\demo{Proof}
Denote $F_s$ the free group generated by $\{ y_1, \dots, y_s \}$, and
the map $r \: F_s \rightarrow F$ sending $y_i$ to $r_i$. The map $d_0 \:
\oplus \R r_i \rightarrow J_F/J_F^2 \cong \gud F \otimes \R$ is the map
induced by $r$, $\gud (r) \otimes \R \: \gud F_s \otimes \R \rightarrow
\gud F \otimes \R$. Furthermore $\ker (\gud (r) \otimes \R) \cong \ker (
\gud (r)) \otimes \R$, as $\gud F_s$ is a free abelian group. Thus $\ker
d_0$ admits a basis $\bar w_1, \dots, \bar w_k$, with the $w_i$ words in
$F_s$ mapping to $\Gamma_2 F$ by $r$.

Now, the map $\Gamma_2 F \rightarrow J_F^2$ sends a bracket $(a,b)$ to
$(a-1)(b-1)-(b-1)(a-1)+ \text{terms in }J_F^3$, and a product $\prod
(a_i,b_i)$ to $\sum (a_i-1)(b_i-1)-(b_i-1)(a_i-1)+ \text{terms in
}J_F^3$. Therefore, all the $w_i= \prod (a_{j_i}, b_{j_i})$ map to $J_F
\wedge J_F+ J_F^3$.
\qed \enddemo

Lemma 1.9 allows us to define a map $d_1 \: \ker d_0 \rightarrow
\bigwedge^2 \gud G \otimes \R$ by composing
$$
\ker d_0 \longrightarrow (J_F \wedge J_F+ J_F^3)/J_F^3 \cong \bigwedge^2
\gud F \otimes \R \longrightarrow \bigwedge^2 \gud G \otimes \R
$$

\proclaim{Proposition 1.10}
$\text{coker}\, d_1 \cong \gdt G \otimes \R$.
\endproclaim

\demo{Proof}
As we have previously explained, $\gud G \otimes \R \cong J_G/J_G^2
\cong J_F/(J_F^2+K) \cong J_F/(J_F^2+D)$. Thus $\bigwedge^2 \gud G
\otimes \R \cong \left( J_F \wedge J_F+ (J_F^3+ J_F \cdot D+ D \cdot
J_F) \right)/(J_F^3+ J_F \cdot D+ D \cdot J_F)$. Also $f( \ker d_0)= D
\cap J_F^2 \subset J_F \wedge J_F+ J_F^3$ by Lemma 1.9, so
$$
\align
\text{coker}\, d_1 &\cong \left( J_F \wedge J_F+ J_F^3+ J_F \cdot D+ D
\cdot J_F+ D \cap J_F^2 \right)/ \left( J_F^3+ J_F \cdot D+ D \cdot J_F+
D \cap J_F^2 \right) \\
&\cong \left( J_F \wedge J_F+ J_F^3+ K \cap J_F^2 \right)/ \left( J_F^3+
K \cap J_F^2 \right) \cong (J_G \wedge J_G + J_G^3)/J_G^3 \cong \gdt G
\otimes \R
\endalign
$$
the last isomorphism being given by Cor. 1.8.
\qed \enddemo

We are now able to determine the structure of the 2-step nilpotent Lie
algebra $\grl G$ of a finitely presented group $G= \langle x_1, \dots,
x_n \; ; \; r_1, \dots, r_s \rangle$:

\proclaim{Proposition 1.11}
Let $\bigwedge^{\le 2} \gud G \otimes \R$ be the free exterior algebra
generated by $\gud G \otimes \R$ modulo the ideal $\bigwedge^{\ge 3} G
\otimes \R$ generated by wedges of length 3 or more. There is an
isomorphism
$$
\grl G \cong (\bigwedge^{\le 2} \gud G \otimes \R) / (\ker d_0 / \ker
d_1 )
$$
\endproclaim

\demo{Proof}
There is an obvious map of exterior algebras, which is a linear
isomorphism in every degree by the above results.
\qed \enddemo

Thus $\grl G$ is the quotient of a free 2-step nilpotent $\R$-Lie
algebra $\bigwedge^{\le 2} H_1(G; \R)$ by a subspace of 2-brackets $\ker
d_0/ \ker d_1$, which corresponds to the relations of the holonomy
algebra. We have
stated in Ex. 1.4 the case of free groups. Let us examine this structure
in some other simple cases:

\proclaim{Corollary 1.12}
Let $G= \langle x_1, \dots, x_n \; ; \; r \rangle$ be a group admitting
a presentation with a single relation. Then:
\item{(i)} If $r \not\in \Gamma_2 F$, there is an isomorphism $\grl G
\cong \grl F_{n-1}$ with $F_{n-1}$ a free group of rank $n-1$.
\item{(ii)} If $r \in \Gamma_2 F \setminus \Gamma_3 F$, there is an
isomorphism $\grl G \cong \grl F / d_1 (r)$.
\item{(iii)} If $r \in \Gamma_3 F$, there is an isomorphism $\grl G
\cong \grl F$.
\endproclaim

\demo{Proof}
All cases are found by applying Prop. 1.11.
\item{(i)} In this case $\gud G \otimes \R \cong \gud F_{n-1} \otimes
\R$, and as $r \not\in \Gamma_2 F$, $\ker d_0= \{ 0 \}$.
\item{(ii)} In this case the map $F \rightarrow G$ induces an
isomorphism $\gud F \otimes \R \cong \gud G \otimes \R$, $\ker d_0=
\R r$, and as $r \not\in \Gamma_3 F$, the coincidence of the lower
central series and augmentation ideal power filtrations in free groups
([MKS], 5.12,[S]) shows that $r-1 \not\in J_F^3$, hence $d_1 (r) \ne 0$.
\item{(iii)} In this case, $\ker d_0 = \R r$ and again by the above
coincidence of filtrations, $d_1 (r)=0$.
\qed \enddemo

\proclaim{Corollary 1.13}
Let $G= \langle x_1, \dots, x_n \; ; \; r_1, \dots, r_s \rangle$ be a
finitely presented group such that its defining relations may be divided
in two sets: $\{ r_1, \dots, r_k \}$ such that $\bar r_1, \dots, \bar
r_k$ are linearly independent in $\gud F \otimes \R$ and $\{ r_{k+1},
\dots, r_s \}$ which belong to $\Gamma_3 F$. Then there is an
isomorphism $\grl G \cong \grl F_{n-k}$, where $F_{n-k}$ is a free group
of rank $n-k$.
\endproclaim

\demo{Proof}
In this case $\gud G \otimes \R$ has rank $n-k$, $\ker d_0= \R r_{k+1}
\oplus \dots \oplus \R r_n$ because those $r_j$ are commutators and the
other relations form a basis of $\text{Im}\,f$, and $\ker d_1= \ker
d_0$ because $r_{k+1}, \dots, r_n \in \Gamma_3 F$.
\qed \enddemo

\demo{Remark 1.14}
We will be interested in this note in which groups $G$ have a free
2-step nilpotent Lie algebra $\grl G$, which by Prop. 1.11 is equivalent
to $\ker d_0= \ker d_1$.

Generic presentations with less relations than generators produce a free
$\grl G$. The reason is that given a group presentation $G= \langle x_1,
\dots, x_n \; ; \; r_1, \dots, r_s \rangle$ with a number of relations
$s \le n$, $\ker d_0=0$ and therefore $\grl G$ is free, unless the
classes $\bar r_1, \dots, \bar r_s \in \gud F_n \otimes \R$ are linearly
dependent. But the sets of linearly dependent $\bar r_1, \dots, \bar
r_s$ form a codimension $n-s+1$ closed subset of $(\gud F_n \otimes
\R)^s$.

The hypotheses of Corollary 1.13
may be weakened by requiring only that $\{ r_1, \dots, r_k \}$ map on a
basis of $\text{Im}\,d_0$, and the remaining relations $\{ r_{k+1},
\dots, r_s \}$ belong to $\Gamma_3 F \dot N_k$, where $N_k$ is the
normal closure in $F$ of $\{ r_1, \dots, r_k \}$.
\enddemo

\demo{Remark 1.15}
Let us conclude this section by bounding the dimension of $\gdt G
\otimes \R$, which will be examined in the coming sections. It follows
from Prop. 1.10 that
$$
\align
\dim \gdt G \otimes \R &= \binom{\dim \gud G \otimes \R}{2} - \dim \ker
d_0 + \dim \ker d_1         \\
&\ge \binom{\dim \gud G \otimes \R}{2}- \dim \ker d_0
\endalign
$$
\enddemo

\vskip 1cm

\flushpar {\bf \S2. Sullivan's 1-minimal models, brackets and cup
products}

\bigpagebreak

We sum up for the reader's convenience some basic facts on Sullivan's
1-minimal models, its equivalence with the Malcev completion of the
fundamental group and its relation with cup products.

\medpagebreak

Let $X$ be now a differentiable manifold, and $\Cal E^*(X)$ its De Rham
complex.

The theory of minimal models developed by Sullivan shows
that $\Cal E^*(X)$ has a 1-minimal model. This is a certain minimal
commutative differential graded algebra (cdga)
$M(2,0)(X)$, defined as the limit of an inductive system of cdga $M(1,1)
\hookrightarrow M(1,2) \hookrightarrow M(1,3) \hookrightarrow \dots$,
together with an algebra morphism $\rho \: M_X \rightarrow \Cal E^*(X)$
such that $H^0(\rho), H^1(\rho)$ are isomorphisms and $H^2(\rho)$ is a
monomorphism. For a more complete and elementary account of the theory
of minimal models we remit ourselves to [GM], which treats
1-minimal models in Chap. XII.
We will only review the construction up to the second step $M(1,2)$,
which will be used to relate $\pi_1(X,*)$ and cup products on $H^1(X)$.

Define $M(1,1)=\bigwedge(V_1^1)$, with $V_1^1 = H^1(X)$. Every element
of $V_1^1$ is defined to have degree one and boundary zero, and the map
$\rho \: M(1,1) \rightarrow \Cal E^*(X)$ sends every $x \in V_1^1$
to its image in a prefixed $\R$-vector space section $H^1(X) \rightarrow
\text{(cocycles)}^1$.

The (1,2)-minimal model is defined as an extension of $M(1,1)$:
$M(1,2)= \bigwedge (V_1^1 \oplus V_2^1)$, where $V_2^1= \ker (H^2 \rho
\: H^2M(1,1) \rightarrow H^2 \Cal E^*(X))$. For any $v \in V_2^1$ we
define $dv$ as the element of $V_1^1 \wedge V_1^1$ defining its
cohomology class, and if $dv= \sum x_iy_i$, $\rho (v)$ is a linearly
varying primitive of $\sum \rho(x_i) \rho (y_i)$ in $\Cal E^*(X)$.

\demo{Remark 2.1} By its definition, $H^2M(1,1) \cong H^1(X) \wedge
H^1(X)$,
and as $\rho$ is a cdga morphism, there is an isomorphism $V_2^1 = (\ker
H^2M(1,1) \rightarrow
H^2 \Cal E^*(X)) \cong \ker (\cup \: H^1(X) \wedge H^1(X) \rightarrow
H^2(X)$.
\enddemo

The following steps $M(1,n)$ are constructed in a likewise manner,
defining $V_n^1$ as $\ker H^2M \allowmathbreak (1,n-1) \rightarrow H^2
\Cal E^*(X)$, and
$d, \rho$ on it as on $V_2^1$. The inductive limit is denoted $M(2,0)$
and is the 1-minimal model of $X$.

\demo{Remark 2.2}
Among the properties of the 1-minimal model let us remark:
\item - It is well defined up to isomorphism.
\item - It is functorial up to homotopy, i.e., any cdga morphism $\Cal
E^{\ast} Y \rightarrow \Cal E^{\ast} X$ may be lifted to a morphism
$M(2,0)(Y) \rightarrow M(2,0)(X)$, and all its liftings are homotopic in
the cdga category.
\item - As a consequence of its uniqueness, if a map $f \: X \rightarrow
Y$ induces isomorphisms $H^0f, H^1f$ and a monomorphism $H^2f$, then it
induces an isomorphism of 1-minimal models $M(2,0)(Y)
\allowmathbreak \cong M(2,0)(X)$.
\item - The 1-minimal model of $X$ may be computed replacing
$\Cal E^*(X)$ by any other cdga quasi-isomorphic to it.
Thus if $X$ is a complex manifold, we may compute it
from its holomorphic De Rham complex, logarithmic
complexes, Dolbeault complexes, etc..
\enddemo

\medpagebreak

We recall now the dualizing process between Lie algebras and
free commutative differential graded algebras
generated by elements of degree one.

Let $L$ be a finite-dimensional $\R$-Lie algebra. Its
bracket is a bilinear alternating map
$$
[.,.] \: L \wedge L \longrightarrow L
\tag 2.3
$$
Dualizing on both sides, $[.,.]$ has an adjoint map
$$
d \: L^{\vee} \longrightarrow L^{\vee} \wedge L^{\vee}
\tag 2.4
$$
The map $d$ may be extended as a graded derivation to the free graded
algebra $\bigwedge L^{\vee}$, defining the degree of elements in
$V^{\vee}$ to be one. The Jacobi identity satisfied by $[.,.]$ dualizes
then as $d^2=0$.

Reciprocally, if $M= \bigwedge W$ is a free cdga and $\text{deg}\,W=1$,
the differential restricts to a map $d \: W=M^1 \rightarrow M^2=W
\wedge W$, which dualizes to a map $[.,.] \: W^{\vee} \wedge W^{\vee}
\rightarrow W^{\vee}$, and the fact $d^2=0$ in $M$ translates as the
Jacobi identity in $W^{\vee}$.

\demo{Definition 2.5} A Lie algebra $L$ and a free cdga generated by
elements of degree one are {\bf dual} when each one yields the other by
the above processes.
\enddemo

The following result is due to D. Sullivan. The reader will find a
complete proof of it in [BG].

\proclaim{Theorem 2.6}(Sullivan)
Let $X$ be an arc-connected topological space with a finitely presented
fundamental
group $\pi_1(X,*)$. The inductive system $M(1,1) \hookrightarrow M(1,2)
\hookrightarrow \dots $ formed by the (1,n)-minimal models of $X$ and
the projective system $\dots \rightarrow \ld \pi_1X \rightarrow \Cal
L_1 \pi_1X$ described in (1.1) are dual.
\endproclaim

This theorem has important consequences for our purposes. The most
obvious is about the duality as vector spaces:

\proclaim{Corollary 2.7}
$V_n^1 \cong (\Gamma_n/\Gamma_{n+1}(\pi_1X) \otimes \R)^{\vee}$
\endproclaim

The duality Lie algebra-cdga also has consequences:

\proclaim{Lemma 2.8}
The diagram
$$
\CD
V_1^1 \wedge V_1^1 @>\cdot>> H^2M(1,2) \\
@V{\rho \wedge \rho}VV        @VV{H^2\rho}V \\
H^1(X) \wedge H^1(X) @>\cup>> H^2(X)
\endCD
$$
is commutative, and the first column is an isomorphism.
\endproclaim

\demo{Remark 2.9}
Lemma 2.8 implies that the cup product between 1-classes is determined
by the brakets in $\grl \pi_1 X$. This is due to the fact that it
factors through $M(1,2)$, wich is dual to $\ld \pi_1 X \hookrightarrow
\grl \pi_1 X$.
\enddemo

Another particular consequence of Theorem 2.6 we will use is:

\proclaim{Corollary 2.10}
$\dim \ker (\cup \: H^1(X) \wedge H^1(X) \rightarrow H^2(X))= \dim \ld
\pi_1X^{(2)}= \dim \gdt \pi_1(X, \ast) \otimes \R$
\endproclaim

\demo{Proof}
Both spaces are isomorphic to $V_2^1$.
\qed \enddemo

\vskip 1cm

\flushpar {\bf \S3. The $\grl$ algebras of K\"{a}hler groups}

\smallpagebreak

In this chapter we proof that groups with a free $\grl$ algebra cannot
be K\"{a}hler (Thm. 3.3), which by Remark 1.14 is the generic situation in
groups with few defining relations.

\medpagebreak

We recall now the definition of the objects of our interest:

\demo{Definition 3.1}
Let $G$ be a group. It is a {\bf K\"{a}hler group} when $G \cong
\pi_1(X,*)$, where $X$ is a compact K\"{a}hler manifold.
\enddemo

Some restrictions on K\"{a}hler groups are well known: They are finitely
presented, their rank is even.
Sullivan's theory of minimal models, completed by Morgan, \dots ,
yields subtle conditions K\"{a}hler groups must satisfy, due to the
formality of compact K\"{a}hler manifolds ([DGMS]), and to the mixed Hodge
structure its 1-minimal model supports, which is an extension of the
Hodge structure on the cohomology ring ([M1],[M2]). The formality
condition it imposes is in our case:

\proclaim{Proposition 3.2}([M1],9.4)
Let $G$ be a K\"{a}hler group. The Lie algebra $\text{Gr}\,
\Cal L G$ is the quotient of a free Lie algebra by an ideal generated by
sums of length two brackets.
\endproclaim

This is equivalent to the Malcev and holonomy algebras of $G$ being
isomorphic (cf. [Ko]). For instance,
the group $G= \langle x,y \; ; \; ((x,y),y) \rangle$ cannot be K\"{a}hler
because $\Cal L G$ is the quotient of a free algebra by the ideal
generated by a length three bracket $[[\bar x, \bar y], \bar y]$.

In this note we will look in another direction through the results of
the preceding
section. Johnson and Rees rule out in [JR] free groups as K\"{a}hler groups
because they have finite index subgroups of odd rank. In this paragraph
we establish this by studying cup products, and in this way the result
holds for groups with a free 2-step nilpotent Lie algebra $\grl$. This
is a consequence of the Lefschetz decomposition and the nondegenerate
alternating pairing $Q$ that the real cohomology $H^{\ast}(X)$ of a
compact K\"{a}hler manifold $X$ supports. Our standard reference for these
properties will be [W], Chap. V.

\proclaim{Theorem 3.3}
Let $G$ be a finitely presented group such that $\grl G \cong \grl F_n$
for some $n$. Then $G$ is not K\"{a}hler.
\endproclaim

\demo{Proof}
Suppose $\pi_1(X,*) \cong G$.
As $\grl G \cong \grl F_n$, and using Example 1.4, $\dim
\gud G \otimes \R= \dim \gud F_n \otimes \R= n$, and $\dim \gdt G
\otimes \R= \dim \ld G^{(2)}= \dim \grl G^{(2)}= \dim \grl F_n =
\frac{n(n-1)}{2}$. Thus by Cor. 2.10 $\ker \cup= H^1(X) \wedge H^1(X)$,
all cup products are zero. But cohomology classes of rank one are
primitive, and the pairing $Q$
is nondegenerate (see [W],5.6). This leads to a contradiction, hence $G$
cannot be K\"{a}hler.
\qed \enddemo

To use Theorem 3.3, one needs to know when $\grl G \cong \grl F_n$.
This can be established from a presentation of $G$ as we have shown in
section \S1 (see Remark 1.14).

\proclaim{Corollary 3.4}
Let $G= \langle x_1, \dots, x_n \; ; \; r \rangle$ be a K\"{a}hler group
presented with a single relation. Then $r \in \Gamma_2 F \{x_1, \dots,
x_n \} \setminus \Gamma_3 F \{ x_1, \dots, x_n \}$.
\endproclaim

\demo{Proof}
By Corollary 1.12, $r \in \Gamma_2 F \setminus \Gamma_3 F$ is the only
case in which $\grl G$ is not free.
\qed \enddemo

\demo{Examples 3.5}
\item{(i)} ([JR]) Free groups cannot be K\"{a}hler.
\item{(ii)} The group $G= \langle x,y \; ; \; ((x,y),y) \rangle$ cannot
be K\"{a}hler because $((x,y),y) \in \Gamma_3 F$ and therefore $\grl G
\cong \grl F_2$. This group was known to be non-K\"{a}hler by [M1],[M2].
\item{(iii)} The group $G= \langle x,y,z,t \; ; \; x^3y^{-4}z^2y, y^2z^2
\rangle$ is not K\"{a}hler because the two defining relations are
linearly
independent in $\gud F \otimes \R \cong \R^4$, and thus by Cor. 1.13
$\grl G \cong \grl F_2$.
\item{(iv)} The group $G = \langle x_1, \dots, x_5 \; ; \; x_1x_2^2x_1,
x_2x_3^2x_2, x_5x_4^2x_5 \rangle$ has also a linearly independent
relation set, and thus $\grl G \cong \grl F_2$, and $G$ cannot be
K\"{a}hler either.
\item{(v)} Compact Riemann surfaces provide examples showing that
one-relator groups with a defining relation $r \in \Gamma_2 F \setminus
\Gamma_3 F$ are possible.
\enddemo

\vskip 1cm

\flushpar {\bf \S4. One- and two-relator K\"{a}hler groups}

\bigpagebreak

In this chapter we give a lower bound for the number of defining
relations that a presentation of a K\"{a}hler group $\pi_1 X$ must have,
determined by the dimension of the Albanese image of $X$ (Prop. 4.5).
We apply this to fully determine the Malcev completions of one- and
two-relator K\"{a}hler groups in Thm. 4.6. Our starting point 4.3 is due
to F. Campana ([C]).

\medpagebreak

The Albanese variety and Albanese map are another feature of compact
complex manifolds which is very useful to study its fundamental group.
Let us briefly recall it.

\proclaim{Proposition 4.1}
Let $X$ be a compact complex manifold, and let $g= \dim
H^0(\Omega_X^1)$. There is a complex torus $Alb(X)$ (the {\bf Albanese
torus}) and an analytic map $\alpha_X \: X \rightarrow Alb(X)$ (the {\bf
Albanese map}) such that $\alpha_X$ induces an isomorphism $H^1(Alb(X);
\Bbb Z)
@>\cong>> H^1(X; \Bbb Z)/_{\text{Torsion}}$. The pair $(Alb(X),
\alpha_X)$ is determined up to
isomorphism by this property. Moreover $\alpha_X(X)$
is a generating set for $Alb(X)$ as an abelian group.
\endproclaim

Let us fix our notation: Let $X$ be compact K\"{a}hler, $\alpha_X \: X
\rightarrow Alb(X)$ its Albanese map, and denote $Y= \alpha_X(X)$ its
image, which may be singular. We consider a desingularization
$\varepsilon \: \tilde Y \rightarrow Y$, and a desingularization $\tilde
X$ of the pullback of $\alpha_X$:
$$
\CD
\tilde X @>\tilde \alpha_X>> \tilde Y \\
@V{\varepsilon_X}VV        @VV{\varepsilon}V  \\
X @>>\alpha_X>Y
\endCD
\tag 4.2
$$

It is clear that the manifold $\tilde X$ is also compact K\"{a}hler, and
that the map $\varepsilon_X$
induces an isomorphism $\varepsilon_{X*} \: \pi_1\tilde X \rightarrow
\pi_1 X$.

We will call the map $\tilde \alpha_X \: \tilde X \rightarrow \tilde Y$
a {\bf smoothing of the Albanese map} of X. The properties of the
original Albanese map $\alpha_X$ relate $X$, $\tilde X$ and $\tilde Y$
for our purposes:

\proclaim{Proposition 4.3} ([C])
Let $X$ be compact K\"{a}hler, and $\tilde \alpha_X \: \tilde X \rightarrow
\tilde Y$ be a smoothing of the Albanese map of $X$. Then $\tilde
\alpha_X$ induces an isomorphism $\Cal L(\pi_1 \tilde X) @>\sim>> \Cal
L(\pi_1 \tilde Y)$.
\endproclaim

\demo{Proof}
(Cf. [C]) As $\varepsilon_X$ induces an isomorphism of fundamental
groups, hence $\varepsilon^{\ast} \: H^1(X) \rightarrow H^1( \tilde
X)$ is also an isomorphism. This implies that $Alb(X)$ is the Albanese
torus of $\tilde X$ and $\alpha_X \circ \varepsilon_X= \varepsilon
\circ \tilde \alpha_X$ its Albanese map. As a consequence $\tilde
\alpha_X^{\ast} \: H^1(\tilde Y) \rightarrow H^1(\tilde X)$ is onto. As
$\alpha_X$ itself is also onto, $\tilde \alpha_X^{\ast}$ is also
injective for $H^{\ast}$.  Therefore $\tilde \alpha_X$ induces an
isomorphism
$H^1(\tilde Y) \cong H^1(\tilde X)$ and an injection $H^2(\tilde Y)
\hookrightarrow H^2(\tilde X)$. Thus by universality of the 1-minimal
model, $\tilde \alpha_X$ induces an isomorphism $M(2,0)(\tilde Y) \cong
M(2,0)(\tilde X)$; our statement is its dual.
\qed \enddemo

Thus the study of Malcev completions of K\"{a}hler groups may be reduced to
the study of smoothings of its Albanese images. But rather than follow
this line, we will derive from it consequences on $H^*(X)$, which
will determine Malcev completions of K\"{a}hler groups with one or two
defining relations.

\proclaim{Lemma 4.4}
Let $X$ be compact K\"{a}hler, $Y$ the Albanese image of $X$
and $m= \dim_{\Bbb C}\tilde Y$. Then the graded algebra $H^*(X; \Bbb C)$
contains a free graded exterior algebra $\bigwedge (V)$, where $V$ is a
complex vector space of dimension $m$ and degree 1 spanned by
holomorphic forms.
\endproclaim

\demo{Proof} (cf. [Be] V.18)
Let $y \in Alb(X)$ be a regular point of the Albanese image $Y=
\alpha_X(X)$. As $\dim Y=m$, there are local coordinates $u_1, \dots
u_n$ of $Alb(X)$ in a neighbourhood $U$ of $y$ such that $Y \cap U$ is
defined as $u_{m+1}=0, \dots, u_n=0$. The holomorphic forms $du_1,
\dots, du_m$ are defined on $U$ and, as $Alb(X)$ is parallelizable, the
forms in $\bigwedge (du_1, \dots, du_m)$ extend to global holomorphic
forms on $Alb(X)$. Its direct image $\bigwedge (\alpha_X^* du_1, \dots,
\alpha_X^* du_m)$ defines a subalgebra of holomorphic cohomology classes
in $H^*(X)$ which is free over $y$, hence is free.
\qed \enddemo

The above Lemma together with the correspondence of 2.10 may be used to
bound from below the number of defining relations for K\"{a}hler groups, and
to study those admitting a one- or two-relation presentation.

\proclaim{Proposition 4.5}
Let $G$ be a K\"{a}hler group, $X$ a compact K\"{a}hler manifold such that
$\pi_1 X \cong G$ and $Y$ its Albanese image. Then:
\item{(i)} If $\dim Y=1$, there is an isomorphism $\Cal L G \cong \Cal
L \pi_1 C_g$ with $C_g$ a compact Riemann surface of genus $g$.
\item{(ii)} If $\dim Y=m >1$, $\dim \ker \left( d_0 \: \R r_1 \oplus
\dots
\oplus \R
r_s \rightarrow \gud F \otimes \R \right) \ge 2\binom{m}{2}+1$. In
particular, any presentation
$G= \langle x_1, \dots, x_n \; ; \; r_1, \dots, r_s \rangle$ must have
defining relations $r_1, \dots, r_k$ such that they form a basis of
$\text{Im}\,f$ and at least another $2 \binom{m}{2}+1$ defining
relations.
\endproclaim

\demo{Proof}
\item{(i)} is just Prop. 4.3 with $\tilde Y$ as $C_g$.
\item{(ii)} By Lemma 4.4, the algebra $H^*(X;\C)$ contains a free
algebra
$\bigwedge(V)$ generated by $m$ linearly independent holomorphic
1-forms.
By the Hodge structure of $H^*(X)$ it contains an isomorphic algebra
$\bigwedge (\bar V)$ spanned by $m$ independent antiholomorphic 1-forms.
Both algebras being free, one obtains the lower bound $\dim \text{Im}\,
\cup \: H^1(X) \wedge H^1(X) \rightarrow H^2(X) \ge 2 \binom{m}{2}$
considering either holomorphic or antiholomorphic products alone.
Finally, due to the properties of the $Q$ pairing in $H^1(X)$ ([W] 5.6),
the product of a holomorphic 1-form with its conjugate cannot be
zero, so $\dim (\text{Im}\,\cup) \cap H^{1,1}(X) \ge 1$.
By the correspondence of 2.10 this produces the sought bound.
\qed \enddemo

We are now able to complete our study of compact K\"{a}hler groups with one
or two defining relations begun in Cor. 3.4.

\proclaim{Theorem 4.6}
Let $G$ be a K\"{a}hler group admitting a presentation with only
one or two defining relations. Then either
$\gud G \otimes \R=0$ or $\Cal L G
\cong \Cal L \pi_1 C_g$ with $C_g$ a compact Riemann surface.
\endproclaim

\demo{Proof}
If $\gud G \otimes \R \ne 0$, then $\grl G \ne 0$, and by Prop. 4.5
any presentation of $G$ must have at least $2\binom{ \dim Y}{2}+1$
defining relations, with $Y= Alb (X)$. Thus the only possible case
is $\dim Y=1$, and Prop. 4.5 (i) completes the proof.
\qed \enddemo

\demo{Remark}
\item - The 1-relator groups $G$ with $\gud G \otimes \R=0$ are
exactly the $G \cong \Bbb Z/ n \Bbb Z$.
\item - The 2-relator groups $G$ with $\gud G \otimes \R=0$ are those
with
a presentation $\langle x_1, x_2 \; ; \; r_1, r_2 \rangle$ with $\bar
r_1, \bar r_2$ linearly independent in $\gud F\{x_1,x_2\}$.

This is immediately derived from the exact sequence (1.5).
\enddemo

\demo{Examples 4.7} Denote $C_g$ a compact Riemann surface of genus $g$.
\item{(i)} The group $G$ defined in Ex. 3.5 (iii) can also be seen not
to be K\"{a}hler by Thm. 4.6, as $\gud G \otimes \R \cong \R^2$ but
$\grl G \not\cong \grl \pi_1(C_1)$.
\item{(ii)} The group $G = \langle x_1, x_2, x_3, x_4 \; ; \;
(x_1x_2,x_3^2), (x_1x_3x_1,x_4^3) \rangle$ cannot be K\"{a}hler because
$\gud G \otimes \R \cong \R^4$ but $\dim \gdt G \otimes \R= 4 \ne 5
=\dim \gdt \pi_1(C_2)$
\enddemo

\vskip 1cm

\flushpar {\bf \S5. Nonfibered K\"{a}hler groups}

\bigpagebreak

Here we establish a dicothomy between fibered and nonfibered K\"{a}hler
groups, arising from a result by A. Beauville and Y.T. Siu on the
existence of irregular pencils on compact K\"{a}hler manifolds. We skip
the fibered case,
and we give in Prop. 5.7 an upper bound for $\dim \gdt G \otimes
\R$ in the case of nonfibered groups. This translates as a lower bound
for the number of relations that their presentations must have.

\medpagebreak

Let $G= \pi_1(X,\ast)$ be a fundamental group. By
Corollary 2.10 $\dim \gdt G \otimes \R= \dim H^1(X) \wedge H^1(X)- \dim
\text{Im}\, \left( \cup \: H^1(X) \wedge H^1(X) \rightarrow H^2(X)
\right)$. We have seen in \S3 that if $X$ is compact K\"{a}hler,
$\text{Im}\, \cup$ must be nonzero. Now we will establish a lower bound
on its dimension in the case of nonfibered manifolds, by recalling
a result of Castelnuovo-De Franchis and its extension to arbitrary
dimension.

\demo{Definition 5.2}
Let $G$ be a K\"{a}hler group.
\item{(i)} We call $G$ a {\bf fibered K\"{a}hler group} when
$G= \pi_1(X,\ast)$ with $X$ compact K\"{a}hler admitting a nonconstant
holomorphic map $f \: X \rightarrow C_g$, with $C_g$ a compact Riemann
surface of genus $g \ge 2$.
\item{(ii)} We call $G$ a {\bf nonfibered K\"{a}hler group} when $G=
\pi_1(X, \ast)$ with $X$ compact K\"{a}hler not admitting any
nonconstant holomorphic map to a compact Riemann surface of genus $g \ge
2$.
\enddemo

A. Beauville and Y.T. Siu independently proved that the above
definitions make sense:

\proclaim{Proposition 5.3} ([Ca] Appendix,[Siu])
Let $X$ be a compact K\"{a}hler manifold, and $G= \pi_1(X,\ast)$. Then
$X$ admits a nonconstant holomorphic map to a compact Riemann surface of
a given
genus $g \ge 2$ if and only if there is an epimorphic group morphism $G
\rightarrow \pi_1(C_g, \ast)$, with $\pi_1 (C_g, \ast)$ the fundamental
group of a compact Riemann surface of genus $g$.
\endproclaim

Prop. 5.3 means that a K\"{a}hler group $G$ is either fibered or nonfibered,
and that the former are characterised by admitting a $\pi_1(C_g)$ as a
quotient.

\demo{Remark}
If we have an onto map $G \rightarrow H \rightarrow 1$, it induces onto
maps $\Gamma_n/ \Gamma_{n+1} G \otimes \R \rightarrow \Gamma_n/
\Gamma_{n+1} H \otimes \R \rightarrow 0$ for all $n$. This together with
the fact that the lower central series quotients of the $\pi_1C_g$ have
all nonzero rank shows
that nilpotent or rationally nilpotent K\"{a}hler groups must be nonfibered.
\enddemo

We now study the cup products of 1-forms in the case of nonfibered
compact K\"{a}hler manifolds. We begin with an extension of a classical
result (see [Ca]):

\proclaim{Proposition 5.4} (Castelnuovo-De Franchis) Let $X$ be a
compact K\"{a}hler manifold. If there exist $\omega_1, \omega_2$ linearly
independent holomorphic 1-forms such that $\omega_1 \wedge \omega_2=0$
then there is a holomorphic map $f \: X \rightarrow C$ with $C$ a curve
of genus $g(C) \ge 2$, such that
$\omega_1, \omega_2$ belong to $\text{Im}\, f^{\ast}$.
\endproclaim

\demo{Remark}
The form equality $\omega_1 \wedge \omega_2=0$ is
equivalent to $\omega_1 \wedge \omega_2$ being exact. This is a result
of Hodge theory, showing that a nonzero holomorphic form over a compact
K\"{a}hler manifold cannot be exact.
\enddemo

The Castelnuovo-De Franchis theorem together with the conic structure of
the set of products in $H^{2,0}(X)$ yield the following corollary (see
[BPV] IV, Prop. 4.2):

\proclaim{Corollary 5.5}
If $X$ is a nonfibered compact K\"{a}hler manifold, then $\dim \text{Im}\,
\cup \: H^{1,0}(X) \wedge H^{1,0}(X) \rightarrow H^{2,0}(X) \ge 2 \dim
H^{1,0}(X)-3$.
\endproclaim

Cor. 5.5 gives a bound for the products of holomorphic 1-forms, and by
conjugation, of antiholomorphic 1-forms. The dimension of products of
holomorphic-antiholomorphic 1-forms has been bounded for compact
complex surfaces in [BPV], IV, Prop. 4.3. We slightly alter their proof
to extend it to compact K\"{a}hler manifolds of arbitrary dimension:

\proclaim{Proposition 5.6}
Let $X$ be a nonfibered compact K\"{a}hler manifold. Then $\dim \left(
\text{Im}\,
\cup \: \right.$ $\left. H^{1,0}(X) \otimes H^{0,1}(X) \rightarrow
H^{1,1}(X) \right) \ge 2 \dim H^{1,0}(X)-1$.
\endproclaim

\demo{Proof}
Denote $n= \dim X \ge 2$, $V = \text{Im}\, \cup \: H^{1,0}(X) \wedge
H^{0,1}(X) \rightarrow H^{1,1}(X)$ and fix $\omega$ a fundamental
K\"{a}hler form on
$X$. We begin by showing that the pairing $\cup \: H^{1,0}(X) \wedge
H^{0,1}(X) \rightarrow V$ becomes injective when we fix a
nonzero $\xi \in H^{1,0}(X)$ or $\bar \eta \in H^{0,1}(X)$.

Suppose there are holomorphic 1-forms $\xi, \eta$ such that $\xi \wedge
\bar \eta= d \alpha$. Then obviously $\xi \wedge \eta \wedge \bar \xi
\wedge \bar \eta= d \alpha'$, and
$$
\int_X \xi \wedge \eta \wedge \bar \xi \wedge \bar \eta \wedge
\omega^{n-2}=0
$$
By the properties of the pairing $Q$ of compact K\"{a}hler manifolds (see
[W] 5.6), this implies that $\xi \wedge \eta=0$, thus by the
Castelnuovo-De Franchis theorem $\xi$ and $\eta$ are linearly dependent.
Take $\xi= a \eta$, with $a \in \Bbb C^{\ast}$. Then $0= \xi \wedge
\bar \eta= a \eta \wedge \bar \eta$. Again by the properties of the
pairing $Q$, this means that $\xi,\eta=0$.

Thus a map may be defined
$$
\Bbb P(H^{1,0}(X)) \times \Bbb P(H^{0,1}(X)) \longrightarrow \Bbb
P(V)
$$
with injective restrictions fixing a point in either factor of the
source. We apply now the following result from [RV]:

\proclaim{Proposition}
Let $\varphi \: \Bbb P^m(\Bbb C) \times \Bbb P^k(\Bbb C) \rightarrow
\Bbb P^l(\Bbb C)$ be a holomorphic mapping, with $l < m+k$. Then
$\varphi$ factors through one of the projections $\Bbb P^m \times \Bbb
P^k \rightarrow \Bbb P^m$, $\Bbb P^m \times \Bbb P^k \rightarrow \Bbb
P^k$.
\endproclaim

In our case, $\cup$ cannot factor through any of the projections because
it is fiberwise injective in both cases, so it holds that $\dim V \ge 2
\dim H^{1,0}(X)-1$ as was wanted.
\qed
\enddemo

We have now all the required pieces to study $\gdt \otimes \R$ of
nonfibered groups. We return to the notations defined in \S1.

\proclaim{Proposition 5.7}
Let $X$ be a nonfibered compact K\"{a}hler manifold with $q= \dim H^1(X)
= \dim \gud \pi_1(X,\ast)) \otimes \R \ne 0$. Then
$$
\dim \gdt \pi_1(X, \ast) \otimes \R \le \frac{2q(2q-1)}{2}-
2(2q-3)-(2q-1)
$$
\endproclaim

\demo{Proof}
We have seen in Cor. 2.10 that $\dim \gdt \pi_1(X, \ast) \otimes \R=
\dim H^1(X) \wedge H^1(X)- \dim \text{Im}\, \left( \cup \: H^1(X) \wedge
H^1(X) \rightarrow H^2(X) \right)= \frac{2q(2q-1)}{2}- \dim
\text{Im}\,\cup$.

We break $H^1(X)$ into its Hodge components. By Cor. 5.5 $\dim
\left( \text{Im}\,H^{1,0}(X) \wedge H^{1,0}(X) \rightarrow \right.$
$\left. H^{2,0}(X) \right) \ge 2q-3$. The same holds by conjugation for
$H^{0,1}(X) \wedge H^{0,1}(X) \rightarrow H^{0,2}(X)$. Prop. 5.6 gives
the inequality
$ \dim \left( H^{1,0}(X) \wedge H^{0,1}(X) \rightarrow H^{1,1}(X)
\right) \ge 2q-1$ and our statement follows from the addition of bounds.
\qed \enddemo

\proclaim{Corollary 5.8}
Let $G= \langle x_1, \dots, x_n \; ; \; r_1, \dots, r_s \rangle$ be a
finite group presentation, such that the images of $r_1, \dots, r_k$
form a basis of $\text{Im}\,d_0 \cong N/N \cap \Gamma_2F \otimes \R
\hookrightarrow \gud F \otimes \R$. If $G$ is a nonfibered K\"{a}hler group,
the total number of relations must satisfy
$$
s \ge k + 2(n-k-3)+(n-k-1)
$$
\endproclaim

\demo{Proof}
Let us note first that $\dim \gud G \otimes \R= \dim \gud F \otimes \R-
\dim N/N\cap \Gamma_2 F \otimes \R= n-k$.

By Prop. 1.10, there is an exact sequence
$$
0 \longrightarrow \ker d_0/ \ker d_1 \longrightarrow \left( \gud G
\otimes \R \right)^{\wedge 2} \longrightarrow \gdt G \otimes \R
\longrightarrow 0
$$
Applying Cor. 5.7 we obtain that
$$
s-k= \dim \ker d_0 \ge \dim \ker d_0/\ker d_1 \ge 2(n-k-3)+(n-k-1)
$$
\qed \enddemo

\demo{Examples 5.9}
\item{(i)} A group $G= \langle x_1, \dots, x_{2q} \; ; \; w_1, \dots,
w_s \rangle$ with $w_1, \dots, w_s \in \Gamma_2 F$ can be nonfibered
K\"{a}hler only if $s \ge 2(2q-3)+(2q-1)$.
\item{(ii)} {\sl Chain link groups} (see [Ro], 3.G) The group $G_{2q}=
\langle x_1, \dots, x_{2q} \; ; \; (x_1,x_2), \dots,
(x_{2q-1},x_{2q}),\allowmathbreak (x_{2q},x_1) \rangle$ is the
fundamental group of
a link of $2q$ circumferences forming a circular chain, for $q \ge 2$.
This group verifies $k= \dim \gud F \otimes \R- \dim \gud G_{2q} \otimes
\R=0$, and $s=2q< 2(2q-3)+(2q-1)$, and therefore $G_{2q}$ cannot be
nonfibered K\"{a}hler. Broadly speaking, if a link is not very
intertwined, its group is not going to be nonfibered K\"{a}hler.
The group $G_4$ verifies that $\dim \gdt G_4 \otimes \R=2$, and
therefore it cannot be fibered K\"{a}hler either, as it cannot map onto
$\pi_1(C_g, \ast)$ for any $g \ge 2$. The groups $G_{2q}$ with $q \ge 3$
do admit onto mappings to $\pi_1(C_2,\ast)$, and the author does not
know if they are fibered K\"{a}hler.
\item{(iii)} The fundamental group of a compact Riemann surface of genus
$g \ge 2$ admits a presentation $\langle a_1, \dots, a_g, b_1, \dots,
b_g \; ; \; (a_1,b_1) \dots (a_g,b_g) \rangle$. In this presentation
$k=0$, and $s=1< 2(2g-3)+(2g-1)$. Therefore, it can only be the
fundamental group of a fibered K\"{a}hler manifold. This is a particular
case of Prop. 5.3.
\item{(iv)} Let $G= \langle x_1, \dots, x_5 \; ; \;
x_1^{2}x_2^{-2}x_4^2, (x_1,x_2), (x_2,x_3), (x_3,x_4), (x_4,x_5)
\rangle$. In this case $n=5$, $k=1$ as $\text{Im}\,d_0= \langle 2 \bar
x_1-2 \bar x_2+ 2 \bar x_4 \rangle$, and $s=5 < k+2(n-k-3)+(n-k-1)=6$.
Therefore $G$ cannot be nonfibered K\"{a}hler. The group $G$ cannot either
map onto $\pi_1(C_g)$, with $C_g$ a smooth projective curve of genus $g
\ge 2$ because $\dim \gdt G \otimes \R=2$, $\dim \gdt \pi_1(C_g)
\otimes \R= \frac{2g(2g-1)}{2}-1 \ge 5$, so we reach the conclusion that
$G$ cannot be K\"{a}hler.
\enddemo

\vskip 1cm

\flushpar {\bf Bibliography}

\bigpagebreak

\item{[Be]} A. Beauville, {\sl Surfaces alg\'{e}briques complexes}
($3^{me}$ Ed.), Ast\'{e}risque {\bf 54}. S.M.F., 1978.

\item{[BG]} A. Bousfield, V. Gugenheim, {\sl On PL De Rham theory and
rational homotopy type}, Memoirs of the AMS {\bf 179}. A.M.S., 1976.

\item{[BPV]} W. Barth, C. Peters, A. Van de Ven, {\sl Compact complex
surfaces}, Ergeb. der Math. 3 F. {\bf 4} Springer-Verlag, 1984.

\item{[C]} F. Campana, {\sl Remarques sur les groupes de K\"{a}hler
nilpotents}, C.R.A.S. Paris {\bf 317} (1993) 777-780; to appear in Ann.
Sci. ENS.

\item{[Ca]} F. Catanese, {\sl Moduli and classification of irregular
Kaehler manifolds (and algebraic varieties) with Albanese general type
fibrations}, Invent. Math. {\bf 104} (1991) 263-289.

\item{[Ch]} K.T. Chen, {\sl Iterated integrals of differential forms and
loop space cohomology}, Ann. Math. {\bf 97} (1973) 217-246.

\item{[DGMS]} P. Deligne, P. Griffiths, J. Morgan, D. Sullivan, {\sl
Real Homotopy Theory of K\"{a}hler Manifolds}, Invent. Math. {\bf 29} (1975)
245-274.

\item{[G]} M. Gromov, {\sl Sur le groupe fondamental d'une vari\'{e}t\'{e}
K\"{a}hl\'{e}rienne}, C.R.A.S. Paris {\bf 308} (1989) 67-70.

\item{[GM]} P. Griffiths, J. Morgan, {\sl Rational Homotopy Theory and
Differential Forms}, Progr. in Math. {\bf 16} Birkh\"{a}user 1981.

\item{[HMR]} P. Hilton, G. Mislin, J. Roitberg, {\sl Localization of
Nilpotent Groups and Spaces}, Math. Studies {\bf 15} North Holland 1975.

\item{[JR]} F. Johnson, E. Rees, {\sl On the fundamental group of a
complex algebraic manifold}, Bull. London Math. Soc. {\bf 19} (1987)
463-466.

\item{[JR2]} F. Johnson, E. Rees, The fundamental groups of
algebraic varieties, in S. Jackowsky et altri (Eds.) {\sl Algebraic
Topology Poznan 1989}, LNM {\bf 1474}, Springer-Verlag 1991, 75-82.

\item{[Ko]} T. Kohno, {\sl S\'{e}rie de Poincar\'{e}-Koszul associ\'{e}e
aux groupes de tresses pures}, Invent. Math. {\bf 82} (1985) 57-75.

\item{[M1]} J. Morgan, {\sl The algebraic topology of smooth algebraic
varieties}, Publ. Math. IHES {\bf 48} (1978) 137-204.

\item{[M2]} J. Morgan, {\sl Hodge theory for the algebraic topology of
smooth algebraic varieties}, Proc. Symp. Pure Math. {\bf 32} AMS 1978,
119-128.

\item{[Q]} D. Quillen, {\sl Rational Homotopy Theory}, Annals of Math.
(2) {\bf 90} (1969) 205-295.

\item {[Q2]} D. Quillen, {\sl On the Associated Graded Ring of a Group
Ring}, J. of Algebra {\bf 10} (1968) 411-418.

\item{[RV]} R. Remmert, A. Van de Ven, {\sl Zur Funktionentheorie
homogener komplexer Mannigfaltigkeiten}, Topology {\bf 2} (1963)
137-157.

\item{[Ro]} D. Rolfsen, {\sl Knots and Links}, Math. Lect. Series {\bf
7} Publish or Perish 1990.

\item{[S]} J.-P. Serre, {\sl Lie algebras and Lie groups}, W.A. Benjamin
1965.

\item{[Siu]} Y.T. Siu, Strong Rigidity for K\"{a}hler Manifolds and the
Construction of Bounded Holomorphic Functions, in R. Howe (Ed.), {\sl
Discrete Groups in Geometry and Analysis}, Progr. in Math. {\bf 67}
Birkh\"{a}user 1987, 124-151.

\item{[St]} J.R. Stallings, Quotients of the powers of the augmentation
ideals in a group ring, in L.P. Neuwirth, {\sl Knots, Groups and
3-Manifolds}, Annals of Math. Studies {\bf 84}, Princeton University
Press 1975, 101-118.

\item {[W0]} R.O. Wells, {\sl Comparison of De Rham and Dolbeault
cohomology
for proper surjective mappings}, Pacific J. of Math. {\bf 53} (1974)
281-300.

\item{[W]} R.O. Wells, {\sl Differential Analysis on Complex Manifolds},
GTM 65, Springer-Verlag 1980.

\enddocument